\documentclass[10pt]{article}
\usepackage{cite}
\usepackage{mathrsfs}
\usepackage{amsfonts}
\usepackage{amsmath}
\usepackage{amsfonts,amssymb,color}
\usepackage{dsfont}
\usepackage{curves}
\usepackage{mathrsfs}
\usepackage{pifont}
\usepackage{amssymb}
\allowdisplaybreaks

\numberwithin{equation}{section}

\date{}

\def\BigRoman{\uppercase\expandafter{\romannumeral\number\count 255 }}
\def\Romannumeral{\afterassignment\BigRoman\count255=}

\begin{document}
\title{The existence of odd-even factors in 1-binding graphs
%\thanks{}
}
\author{\small Sizhong Zhou$^{1}$\footnote{Corresponding author. E-mail address: zsz\_cumt@163.com (S. Zhou)}, Qiuxiang Bian$^{1}$, Hongxia Liu$^{2}$\\
\small $1$. School of Science, Jiangsu University of Science and Technology,\\
\small Zhenjiang, Jiangsu 212100, China\\
\small $2$. School of Mathematics and Information Science, Yantai University,\\
\small Yantai, Shandong 264005, China\\
}

\maketitle
\begin{abstract}
\noindent Let $G$ be a graph. The binding number of $G$, denoted by $\mbox{bind}(G)$, is defined as
$$
\mbox{bind}(G)=\min\left\{\frac{|N_G(S)|}{|S|}:\emptyset\neq S\subseteq V(G) \ \mbox{and} \ N_G(S)\neq V(G)\right\}.
$$
If $\mbox{bind}(G)\geq r$, then $G$ is called $r$-binding, where $r$ is a positive real number. The adjacency matrix of $G$ is denoted by $A(G)$. The largest eigenvalue of $A(G)$, denoted by $\rho(G)$,
is said to be the spectral radius of $G$. A spanning subgraph $F$ of $G$ is called an odd-even factor $F=F_W$ if $d_F(u)\in\{1,3,\ldots,k\}$ for every $u\in W$ and $d_F(v)\in\{0,2,\ldots,k+1\}$ for
every $v\in V(G)-W$, where $k$ is a positive odd integer and $W$ is any set of even number of vertices of $G$. In this paper, we propose a tight sufficient condition based on the spectral radius to
guarantee that a connected 1-binding graph $G$ contains an odd-even factor $F=F_W$ such that $d_F(u)\in\{1,3,\ldots,k\} \ \mbox{for all} \ u\in W$ and $d_F(v)\in\{0,2,\ldots,k+1\} \ \mbox{for all} \ v\in V(G)-W$.
\\
\begin{flushleft}

(2020) Mathematics Subject Classification: 05C50, 05C70
\end{flushleft}
\end{abstract}

\section{Introduction}

Throughout this paper, we only deal with simple and undirected graphs. Let $G$ be a graph with vertex set $V(G)$ and edge set $E(G)$. The number of vertices of $G$ is called the order of $G$,
denoted by $|V(G)|=n$. We denote by $\omega(G)$ the number of components in $G$. Given $v\in V(G)$, the set of vertices adjacent to $v$ in $G$ is called the neighborhood of $v$, denoted by $N_G(v)$,
and $d_G(v)=|N_G(v)|$ is said to be the degree of $v$ in $G$. Given a subset $S$ of $V(G)$, we use $G[S]$ and $G-S$ to denote the subgraphs of $G$ induced by $S$ and $V(G)-S$, respectively. Anderson
\cite{As} and Woodall \cite{W} introduced a parameter, binding number of $G$, denoted by $\mbox{bind}(G)$, which is defined as
$$
\mbox{bind}(G)=\min\left\{\frac{|N_G(S)|}{|S|}:\emptyset\neq S\subseteq V(G) \ \mbox{and} \ N_G(S)\neq V(G)\right\}.
$$
If $\mbox{bind}(G)\geq r$, then $G$ is called $r$-binding, where $r$ is a positive real number. We denote by $K_n$ the complete graph with $n$ vertices. Given two disjoint graphs $G_1$ and $G_2$,
the join and the union of $G_1$ and $G_2$ are denoted by $G_1\vee G_2$ and $G_1\cup G_2$, respectively.

The adjacency matrix of $G$ is denoted by $A(G)$. The largest eigenvalue of $A(G)$, denoted by $\rho(G)$, is said to be the spectral radius of $G$. We refer the reader to \cite{GA,B,F,ZZZL,ZZL,Ws2} for
some properties of the spectral radius in a graph.

A spanning subgraph $F$ of $G$ is called an odd-even factor $F=F_W$ if $d_F(u)\in\{1,3,\ldots,k\}$ for every $u\in W$ and $d_F(v)\in\{0,2,\ldots,k+1\}$ for every $v\in V(G)-W$, where $k$ is a
positive odd integer and $W$ is any set of even number of vertices of $G$. Obviously, an odd-even factor is a generalization of an odd (resp. even) factor. Lots of scholars investigated the
connection between spectral radii and factors in graphs. O \cite{Os} established a tight spectral condition for a connected graph with a perfect matching. Fan and Lin \cite{FL} obtained a spectral
condition for a connected 1-binding graph to contain a perfect matching. Zhou, Bian and Sun \cite{ZBS}, Li and Miao \cite{LM} proposed some sufficient conditions based on the spectral radius for
the existence of $[1,2]$-factors in connected graphs. Fan, Goryainov, Huang and Lin \cite{FGHL}, Wu \cite{Wc}, Zhou \cite{Za}, Wu, Zhou and Liu \cite{WZL} provided some spectral radius conditions
for connected graphs to possess $[1,b]$-factors. Zhou and Liu \cite{ZL} characterized a connected graph with an odd $[1,b]$-factor by using the spectral radius. Zhou \cite{Zs},  Wang, Yang and Yang
\cite{WYY} showed some sufficient spectral conditions for the existence of odd $[1,b]$-factors with given properties in connected graphs. Hua and Zhang \cite{HZ} proposed a sufficient condition by
virtue of the distance spectral radius which ensures the existence of an odd $[1,b]$-factor in a connected 1-binding graph. Zhou, Bian and Wu \cite{ZBW} created a spectral radius condition to
guarantee that a connected graph $G$ of even order with minimum degree $\delta$ has an even factor. Li, Liu, Hua and Du \cite{LLHD} put forward a sufficient condition in terms of the signless
Laplacian spectral radius to ensure that a connected graph $G$ of even order with minimum degree $\delta$ has an even factor. We can refer the reader to \cite{Ws1,ZWH,Zt,Zr,Zs,PZ,MLW} for more results
on graph factors.

Lu and Kano \cite{LK} gave a characterization for a graph with an odd-even factor, which is generalized by Kano and Matsumura \cite{KM}. Motivated by \cite{KM,LK}, we study the connection between
spectral radius and odd-even factors in connected 1-binding graphs, and obtain the following result.

\medskip

\noindent{\textbf{Theorem 1.1.}} Let $k$ be a positive odd integer, and let $G$ be a connected 1-binding graph of order $n\geq n(k)$, where
\[
n(k)=\left\{
\begin{array}{ll}
10,&\mbox{if} \ k=1;\\
2k+4,&\mbox{if} \ k\geq3.\\
\end{array}
\right.
\]
If
$$
\rho(G)\geq\rho(K_1\vee(K_{n-2k-2}\cup kK_2\cup K_1)),
$$
then for any set $W$ of even number of vertices of $G$, $G$ contains an odd-even factor $F=F_W$ such that $d_F(u)\in\{1,3,\ldots,k\} \ \mbox{for all} \ u\in W$ and $d_F(v)\in\{0,2,\ldots,k+1\} \ \mbox{for all} \ v\in V(G)-W$,
unless $G=K_1\vee(K_{n-2k-2}\cup kK_2\cup K_1)$.

\section{Preliminary lemmas}

In this section, we introduce some preliminary lemmas, which are used to prove our main theorems. Lu and Kano \cite{LK} established a characterization for a graph having an odd-even factor,
which is generalized by Kano and Matsumura \cite{KM}.

\medskip

\noindent{\textbf{Lemma 2.1}} (Kano and Matsumura \cite{KM}, Lu and Kano \cite{LK}). Let $f:V(G)\rightarrow$\{1,3,5,\ldots\} be a function. Then for any set $W$ of even number of vertices of $G$,
$G$ contains an odd-even factor $F=F_W$ such that
$$
d_F(u)\in\{1,3,\ldots,f(u)\} \ \mbox{for all} \ u\in W, \ \mbox{and}
$$
$$
d_F(v)\in\{0,2,\ldots,f(u)+1\} \ \mbox{for all} \ v\in V(G)-W,
$$
if and only if
$$
\omega(G-S)\leq f(S)+1
$$
for all $S\subset V(G)$.

\medskip

We immediately obtain the following lemma if $f(x)=k$ for any $x\in V(G)$ in Lemma 2.1, where $k$ is a positive odd integer.

\medskip

\noindent{\textbf{Lemma 2.2.}} Let $k$ be a positive odd integer. Then for any set $W$ of even number of vertices of $G$, $G$ contains an odd-even factor $F=F_W$ such that
$$
d_F(u)\in\{1,3,\ldots,k\} \ \mbox{for all} \ u\in W, \ \mbox{and}
$$
$$
d_F(v)\in\{0,2,\ldots,k+1\} \ \mbox{for all} \ v\in V(G)-W,
$$
if and only if
$$
\omega(G-S)\leq k|S|+1
$$
for all $S\subset V(G)$.

\medskip

\noindent{\textbf{Lemma 2.3}} (Li and Feng \cite{LF}). Let $H$ be a subgraph of a connected graph $G$. Then
$$
\rho(G)\geq\rho(H),
$$
with equality occurring if and only if $G=H$.

\medskip

\noindent{\textbf{Lemma 2.4.}} Let $q$ and $s$ be two positive integers with $q\geq s+2$, and let $n_1\geq n_2\geq\cdots\geq n_q\geq1$ be integers with $\sum\limits_{i=1}^{q}n_i=n-s$. If
$n_i\geq2$ for $1\leq i\leq q-s$, $n_j\geq1$ for $q-s+1\leq j\leq q$ and $n_1<n-2q+2$, then
$$
\rho(K_s\vee(K_{n_1}\cup K_{n_2}\cup\cdots\cup K_{n_q}))<\rho(K_s\vee(K_{n-2q+2}\cup(q-s-1)K_2\cup sK_1)).
$$

\medskip

\noindent{\it Proof.} Let $G=K_s\vee(K_{n_1}\cup K_{n_2}\cup\cdots\cup K_{n_q})$, and let $x$ denote the Perron vector of $A(G)$. Using symmetry, we can assume that $x(v)=x_i$ for each
$v\in V(K_{n_i})$, where $1\leq i\leq q$, and $x(u)=y$ for each $u\in V(K_s)$. Since $K_{n_1+s}$ is a proper subgraph of $G$, it follows from Lemma 2.3 and $s\geq1$ that
$$
\rho(G)>\rho(K_{n_1+s})=n_1+s-1>n_1-1.
$$
Notice that $n_1\geq n_i$ for $2\leq i\leq q$. Then from $A(G)x=\rho(G)x$, we possess
$$
(\rho(G)-(n_i-1))(x_1-x_i)=(n_1-n_i)x_1\geq0,
$$
where $2\leq i\leq q$. This leads to $x_1\geq x_i$ for $2\leq i\leq q$. Let $G'=K_s\vee(K_{n-2q+2}\cup(q-s-1)K_2\cup sK_1)$. By means of $n_i\geq2$ for $1\leq i\leq q-s$, $n_j\geq1$ for
$q-s+1\leq j\leq q$, $n_1<n-2q+2$ and $x_1\geq x_i$ for $2\leq i\leq q$, we deduce
\begin{align*}
\rho(G')-\rho(G)\geq&x^{T}(A(G')-A(G))x\\
=&2\sum\limits_{i=2}^{q-s}n_1(n_i-2)x_1x_i+2\sum\limits_{i=q-s+1}^{q}n_1(n_i-1)x_1x_i\\
&-4\sum\limits_{i=2}^{q-s}(n_i-2)x_i^{2}-2\sum\limits_{i=q-s+1}^{q}(n_i-1)x_i^{2}\\
&+2\sum\limits_{i=2}^{q-s-1}\sum\limits_{j=i+1}^{q-s}(n_i-2)(n_j-2)x_ix_j\\
&+2\sum\limits_{i=2}^{q-s}\sum\limits_{j=q-s+1}^{q}(n_i-2)(n_j-1)x_ix_j\\
&+2\sum\limits_{i=q-s+1}^{q-1}\sum\limits_{j=i+1}^{q}(n_i-1)(n_j-1)x_ix_j\\
=&2\sum\limits_{i=2}^{q-s}(n_i-2)(n_1x_1-2x_i)x_i+2\sum\limits_{i=q-s+1}^{q}(n_i-1)(n_1x_1-x_i)x_i\\
&+2\sum\limits_{i=2}^{q-s-1}\sum\limits_{j=i+1}^{q-s}(n_i-2)(n_j-2)x_ix_j\\
&+2\sum\limits_{i=2}^{q-s}\sum\limits_{j=q-s+1}^{q}(n_i-2)(n_j-1)x_ix_j\\
&+2\sum\limits_{i=q-s+1}^{q-1}\sum\limits_{j=i+1}^{q}(n_i-1)(n_j-1)x_ix_j\\
>&0,
\end{align*}
which implies $\rho(G)<\rho(G')$. That is to say,
$$
\rho(K_s\vee(K_{n_1}\cup K_{n_2}\cup\cdots\cup K_{n_q}))<\rho(K_s\vee(K_{n-2q+2}\cup(q-s-1)K_2\cup sK_1)).
$$
Lemma 2.4 is proved. \hfill $\Box$

\medskip

Let $M$ be a real $n\times n$ matrix and $\mathcal{N}=\{1,2,\ldots,n\}$. Assume that the matrix $M$, based on the partition $\pi:\mathcal{N}=\mathcal{N}_1\cup\mathcal{N}_2\cup\cdots\cup\mathcal{N}_r$,
can be written as
\begin{align*}
M=\left(
  \begin{array}{cccc}
    M_{11} & M_{12} & \cdots & M_{1r}\\
    M_{21} & M_{22} & \cdots & M_{2r}\\
    \vdots & \vdots & \ddots & \vdots\\
    M_{r1} & M_{r2} & \cdots & M_{rr}\\
  \end{array}
\right).
\end{align*}
The quotient matrix of $M$ based on $\pi$ is defined by the $r\times r$ matrix $M_{\pi}=(m_{ij})$, where $m_{ij}$ is the average row sum of $M_{ij}$. The partition $\pi$ is called equitable
if the row sum of each block $M_{ij}$ of $M$ is a constant.

\medskip

\noindent{\textbf{Lemma 2.5}} (You, Yang, So and Xi \cite{YYSX}). Let $M$ be a real $n\times n$ matrix with an equitable partition $\pi$, and let $M_{\pi}$ be the corresponding quotient matrix.
Then the eigenvalues of $M_{\pi}$ are also eigenvalues of $M$. Furthermore, if $M$ is nonnegative and irreducible, then the largest eigenvalues of $M$ and $M_{\pi}$ are equal.

\medskip

\noindent{\textbf{Lemma 2.6.}} Let $k$ and $n$ be positive integers with $k\equiv1$ (mod 2) and $n\geq2k+4$. Then $K_1\vee(K_{n-2k-2}\cup kK_2\cup K_1)$ has no odd-even factor $F=F_W$ such that
$d_F(u)\in\{1,3,\ldots,k\}$ for all $u\in W$ and $d_F(v)\in\{0,2,\ldots,k+1\}$ for all $v\in V(K_1\vee(K_{n-2k-2}\cup kK_2\cup K_1))-W$, where $W$ is any set of even number of vertices of
$K_1\vee(K_{n-2k-2}\cup kK_2\cup K_1)$.

\medskip

\noindent{\it Proof.} Let $G=K_1\vee(K_{n-2k-2}\cup kK_2\cup K_1)$. Choose $S=V(K_1)$ in $G$ such that $G-S=K_{n-2k-2}\cup kK_2\cup K_1$. Thus, we possess
$$
\omega(G-S)=k+2=k|S|+2>k|S|+1.
$$
By virtue of Lemma 2.2, $G$ has no desired factor. Lemma 2.6 is verified. \hfill $\Box$

\medskip

\section{The proof of Theorem 1.1}

\noindent{\it Proof of Theorem 1.1.} Suppose that a connected 1-binding $G$ contains no desired factor. In terms of Lemma 2.2, there exists a nonempty subset $S\subset V(G)$ satisfying
$$
\omega(G-S)\geq k|S|+2.
$$
Set $|S|=s$ and $\omega(G-S)=q$. Then $q\geq ks+2$. We denote the components of $G-S$ by $C_1,C_2,\ldots,C_q$, where $|C_i|=n_i$ for $i=1,2,\ldots,q$. Without loss of generality, we let
$n_q\geq n_{q-1}\geq\cdots\geq n_1$. We are to prove $n_{s+1}\geq2$. Otherwise, we obtain $n_i=1$ for $i=1,2,\ldots,s+1$. Set $X=V(C_1\cup C_2\cup\cdots\cup C_{s+1})$. We easily see
$N_G(X)\subseteq S$ and $|X|=s+1$. Thus, we possess
$$
\frac{|N_G(X)|}{|X|}\leq\frac{|S|}{|X|}=\frac{s}{s+1}<1,
$$
which contradicts that $G$ is 1-binding. This leads to $n_i\geq2$ for $i\geq s+1$. Obviously, $G$ is a spanning subgraph of $G_1=K_s\vee(K_{n_1}\cup K_{n_2}\cup\cdots\cup K_{n_{ks+2}})$,
where $n_{ks+2}\geq n_{ks+1}\geq\cdots\geq n_{s+1}\geq2$, $n_s\geq\cdots n_1\geq1$ and $\sum\limits_{i=1}^{ks+2}n_i=n-s$. By virtue of Lemma 2.3, we get
\begin{align}\label{eq:3.1}
\rho(G)\leq\rho(G_1),
\end{align}
with equality if and only if $G=G_1$. Let $G_2=K_s\vee(K_{n-2ks-2}\cup((k-1)s+1)K_2\cup sK_1)$, where $n\geq2ks+4$. Applying Lemma 2.4, we infer
\begin{align}\label{eq:3.2}
\rho(G_1)\leq\rho(G_2),
\end{align}
with equality if and only if $G_1=G_2$.

If $s=1$, then $G_2=K_1\vee(K_{n-2k-2}\cup kK_2\cup K_1)$. Based on \eqref{eq:3.1} and \eqref{eq:3.2}, we obtain
$$
\rho(G)\leq\rho(K_1\vee(K_{n-2k-2}\cup kK_2\cup K_1)),
$$
where the equality occurs if and only if $G=K_1\vee(K_{n-2k-2}\cup kK_2\cup K_1)$, which contradicts the condition of the theorem because $K_1\vee(K_{n-2k-2}\cup kK_2\cup K_1)$ contains no
desired factor (see Lemma 2.6). Next, we shall discuss $s\geq2$.

Recall that $G_2=K_s\vee(K_{n-2ks-2}\cup((k-1)s+1)K_2\cup sK_1)$. The quotient matrix of $A(G_2)$ based on the partition $V(G_2)=V(K_s)\cup V(K_{n-2ks-2})\cup V(((k-1)s+1)K_2)\cup V(sK_1)$
can be written as
\begin{align*}
B_2=\left(
  \begin{array}{cccc}
  s-1 & n-2ks-2 & 2(k-1)s+2 & s\\
  s & n-2ks-3 & 0 & 0\\
  s & 0 & (k-1)s & 0\\
  s & 0 & 0 & 0\\
  \end{array}
\right).
\end{align*}
Based on a simple calculation, the characteristic polynomial of $B_2$ is
\begin{align*}
\varphi_{B_2}(x)=&x^{4}-(n-ks-4)x^{3}+((k-1)sn-n-(2k^{2}-k)s^{2}-(2k-1)s+3)x^{2}\\
&+((2k-1)s^{2}n+(k+1)sn-(4k^{2}-3k+1)s^{3}-(2k^{2}+7k-2)s^{2}-(3k+3)s)x\\
&+(k-1)s^{3}(-n+2ks+3).
\end{align*}
By Lemma 2.5 and the equitable partition $V(G_2)=V(K_s)\cup V(K_{n-2ks-2})\cup V(((k-1)s+1)K_2)\cup V(sK_1)$, we know that $\rho(G_2)$ is the largest root of $\varphi_{B_2}(x)=0$.

Let $G_*=K_1\vee(K_{n-2k-2}\cup kK_2\cup K_1)$. In view of the partition $V(G_*)=V(K_1)\cup V(K_{n-2k-2})\cup V(kK_2)\cup V(K_1)$, the quotient matrix of $A(G_*)$ equals
\begin{align*}
B_*=\left(
  \begin{array}{cccc}
  0 & n-2k-2 & 2k & 1\\
  1 & n-2k-3 & 0 & 0\\
  1 & 0 & k-1 & 0\\
  1 & 0 & 0 & 0\\
  \end{array}
\right).
\end{align*}
Then the characteristic polynomial of $B_*$ is
\begin{align*}
\varphi_{B_*}(x)=&x^{4}-(n-k-4)x^{3}+((k-2)n-2k^{2}-k+4)x^{2}\\
&+(3kn-6k^{2}-7k-2)x+(k-1)(-n+2k+3).
\end{align*}
Applying Lemma 2.5 and the equitable partition $V(G_*)=V(K_1)\cup V(K_{n-2k-2})\cup V(kK_2)\cup V(K_1)$, $\rho(G_*)$ is the largest root of $\varphi_{B_*}(x)=0$.

By a direct computation, we obtain
\begin{align}\label{eq:3.3}
\varphi_{B_2}(x)-\varphi_{B_*}(x)=(s-1)\beta(x),
\end{align}
where $\beta(x)=kx^{3}+((k-1)n-(2k^{2}-k)s-2k^{2}-k+1)x^{2}+((2k-1)sn+3kn-(4k^{2}-3k+1)s^{2}-(6k^{2}+4k-1)s-6k^{2}-7k-2)x-(k-1)n(s^{2}+s+1)+(2k^{2}-2k)s^{3}+(2k^{2}+k-3)s^{2}+(2k^{2}+k-3)s+2k^{2}+k-3$.
We take the derivative function of $\beta(x)$. Notice that $n\geq2ks+4$. For $x\geq n-2k-2$, we deduce
\begin{align*}
\beta'(x)=&3kx^{2}+2((k-1)n-(2k^{2}-k)s-2k^{2}-k+1)x+(2k-1)sn\\
&+3kn-(4k^{2}-3k+1)s^{2}-(6k^{2}+4k-1)s-6k^{2}-7k-2\\
\geq&3k(n-2k-2)^{2}+2((k-1)n-(2k^{2}-k)s-2k^{2}-k+1)(n-2k-2)\\
&+(2k-1)sn+3kn-(4k^{2}-3k+1)s^{2}-(6k^{2}+4k-1)s\\
&-6k^{2}-7k-2 \ \ \ \ \ \ \ \ \ (\mbox{since} \ x\geq n-2k-2)\\
=&(5k-2)n^{2}-((4k^{2}-4k+1)s+20k^{2}+11k-6)n-(4k^{2}-3k+1)s^{2}\\
&+(8k^{3}-2k^{2}-8k+1)s+20k^{3}+30k^{2}+5k-6\\
\geq&(5k-2)(2ks+4)^{2}-((4k^{2}-4k+1)s+20k^{2}+11k-6)(2ks+4)\\
&-(4k^{2}-3k+1)s^{2}+(8k^{3}-2k^{2}-8k+1)s\\
&+20k^{3}+30k^{2}+5k-6 \ \ \ \ \ \ \ \ \ (\mbox{since} \ n\geq2ks+4)\\
=&(12k^{3}-4k^{2}+k-1)s^{2}-(32k^{3}-40k^{2}+12k+3)s\\
&+20k^{3}-50k^{2}+41k-14\\
\geq&4(12k^{3}-4k^{2}+k-1)-2(32k^{3}-40k^{2}+12k+3)\\
&+20k^{3}-50k^{2}+41k-14 \ \ \ \ \ \ \ \ \ (\mbox{since} \ s\geq2)\\
=&4k^{3}+14k^{2}+21k-24\\
>&0 \ \ \ \ \ \ \ \ \ (\mbox{since} \ k\geq1).
\end{align*}
Hence, we claim that $\beta(x)$ is increasing for $x\geq n-2k-2$, and obtain
\begin{align}\label{eq:3.4}
\beta(x)\geq&\beta(n-2k-2)\nonumber\\
=&(2k-1)n^{3}-((2k^{2}-3k+1)s+12k^{2}+4k-5)n^{2}\nonumber\\
&-((4k^{2}-2k)s^{2}-(8k^{3}-6k^{2}-11k+4)s-24k^{3}-28k^{2}+6k+9)n\nonumber\\
&+(2k^{2}-2k)s^{3}+(8k^{3}+4k^{2}-3k-1)s^{2}\nonumber\\
&-(8k^{4}-22k^{2}-11k+5)s-16k^{4}-32k^{3}-8k^{2}+15k+5\nonumber\\
\triangleq&\gamma(n).
\end{align}
By applying a similar analysis as above, we also infer that $\gamma(n)$ is increasing for $n\geq2ks+4$. Thus, we possess
\begin{align}\label{eq:3.5}
\gamma(n)\geq&\gamma(2ks+4)\nonumber\\
=&(2k-1)(2ks+4)^{3}-((2k^{2}-3k+1)s+12k^{2}+4k-5)(2ks+4)^{2}\nonumber\\
&-((4k^{2}-2k)s^{2}-(8k^{3}-6k^{2}-11k+4)s-24k^{3}-28k^{2}+6k+9)(2ks+4)\nonumber\\
&+(2k^{2}-2k)s^{3}+(8k^{3}+4k^{2}-3k-1)s^{2}\nonumber\\
&-(8k^{4}-22k^{2}-11k+5)s-16k^{4}-32k^{3}-8k^{2}+15k+5\nonumber\\
=&(8k^{4}-4k^{3}+2k^{2}-2k)s^{3}-(32k^{4}-44k^{3}+14k^{2}+3k+1)s^{2}\nonumber\\
&+(40k^{4}-104k^{3}+82k^{2}-19k-5)s-16k^{4}+64k^{3}-88k^{2}+55k-15\nonumber\\
\triangleq&p(s).
\end{align}
By utilizing a similar analysis as above, we also ensure that $p(s)$ is increasing for $s\geq2$. For $s\geq3$ and $n\geq2ks+4$, it follows from \eqref{eq:3.5} and $k\geq1$ that
$$
\gamma(n)\geq p(s)\geq p(3)=32k^{4}+40k^{3}+86k^{2}-83k-39>0.
$$
For $s=2$ and $n\geq2ks+4$, it follows from \eqref{eq:3.5} and $k\geq3$ that
$$
\gamma(n)\geq p(s)=p(2)=36k^{2}-11k-29>0.
$$
For $s=2$, $n\geq10$ and $k=1$, it follows from \eqref{eq:3.5} that $\gamma(n)\geq\gamma(10)=126>0$. Combining these with \eqref{eq:3.3}, \eqref{eq:3.4} and $s\geq2$, we obtain
$$
\varphi_{B_2}(x)-\varphi_{B_*}(x)=(s-1)\beta(x)\geq(s-1)\gamma(n)>0
$$
for $x\geq n-2k-2$, which yields
$$
\varphi_{B_2}(x)>\varphi_{B_*}(x)
$$
for $x\geq n-2k-2$. Since $K_{n-2k-1}$ is a proper subgraph of $K_1\vee(K_{n-2k-2}\cup kK_2\cup K_1)$, we deduce $\rho(K_1\vee(K_{n-2k-2}\cup kK_2\cup K_1))>\rho(K_{n-2k-1})=n-2k-2$.
Recall that $\varphi_{B_2}(\rho(G_2))=0$ and $\varphi_{B_*}(\rho(G_*))=0$. Hence, we have
\begin{align}\label{eq:3.6}
\rho(G_2)<\rho(G_*)
\end{align}
for $s\geq2$. It follows from \eqref{eq:3.1}, \eqref{eq:3.2} and \eqref{eq:3.6} that
$$
\rho(G)\leq\rho(G_1)\leq\rho(G_2)<\rho(G_*)=\rho(K_1\vee(K_{n-2k-2}\cup kK_2\cup K_1))
$$
for $s\geq2$, which contradicts $\rho(G)\geq\rho(K_1\vee(K_{n-2k-2}\cup kK_2\cup K_1))$. This completes the proof of Theorem 1.1. \hfill $\Box$

\section*{Declaration of competing interest}

The authors declare that they have no known competing financial interests or personal relationships that could have appeared to influence the work reported in this paper.

\section*{Data availability}

No data was used for the research described in the article.

\section*{Acknowledgments}

This work was supported by the Natural Science Foundation of Jiangsu Province (Grant No. BK20241949). Project ZR2023MA078 supported by Shandong Provincial Natural Science Foundation.

\end{document}